\documentclass{amsart}
\usepackage{amssymb}
\usepackage{amsmath,amscd}
\usepackage{combelow}

\newcommand{\X}{\mathfrak{X}}
\newcommand{\arr}{\rightarrow}
\newcommand{\pb}{\{ \cdot, \cdot \}}

\newtheorem{theorem}{Theorem}

\newtheorem{claim}[theorem]{Claim}

\newtheorem*{lemma*}{Lemma}
\newtheorem*{theorem*}{Theorem}
\newtheorem*{remark*}{Remark}

\theoremstyle{remark}

\theoremstyle{definition}

\newtheorem*{definition*}{Definition}
\newtheorem*{claim*}{Claim}

\begin{document}

\title[Non-degeneracy of the Hofer norm for Poisson structures]{Non-degeneracy of the Hofer norm for Poisson structures}
\author{ Du\v{s}an Joksimovi\'c, Ioan M\u{a}rcu\cb{t}}
\address{Sorbonne Universit\'e and Universit\'e de Paris, CNRS, IMJ-PRG}
\email{joksimovic@imj-prg.fr}
\address{Radboud University Nijmegen}
\email{i.marcut@math.ru.nl}

\begin{abstract}
We remark that, as in the symplectic case, the Hofer norm on the Hamiltonian group of a Poisson manifold is non-degenerate. 
The proof follows from the symplectic case after reducing the problem to a symplectic leaf.
\end{abstract}

\maketitle

A Poisson structure on a smooth manifold $M$ is given by a Lie bracket $\pb$ on the space $C^{\infty}(M)$ which satisfies the Leibniz identity
$$\{fg,h\} = f\{g,h\}+ g\{f,h\}, \quad \forall f,g,h \in C^{\infty}(M).$$
A smooth map between two Poisson manifolds $\psi: (N,\pb_N) \arr (M, \pb_M)$ is called Poisson if the following holds
$$\{f,g\}_M \circ \psi = \{f \circ \psi, g \circ \psi\}_N, \quad \forall f,g \in C^{\infty}(M).$$
Examples of Poisson manifolds include symplectic manifolds, duals of Lie algebras, and every manifold carries the trivial Poisson structure $\{\cdot,\cdot \} \equiv 0.$ 

Let $(M,\pb)$ be a Poisson manifold and $f \in C_c^{\infty}([0,1] \times M)$ be a compactly supported time-dependent Hamiltonian function. We define the (time-dependent) Hamiltonian vector field $X_f^t$ associated to $f$ by
$$X_f^t := \{f_t, \cdot\} \in \X (M),$$
where $f_t := f(t, \cdot) \in C^{\infty}(M),$ $t\in [0,1].$ 
The flow $\{\varphi_f^t\}$  of $X_f^t$ is called the Hamiltonian flow (or the Hamiltonian isotopy) generated by $f.$ 
The Hamiltonian group of $(M,\pb)$ is
\[\mathrm{Ham}(M,\pb):=\big\{\varphi^1_f\ : \ f\in C^{\infty}_{c}([0,1]\times M)\big\}.\]
The length of the Hamiltonian isotopy generated by $f$ is defined as
\[l(f):=\int_0^1\big(\sup_{x \in M} f_t(x)-\inf_{x \in M} f_t(x) \big) \ dt.\]
Notice that, in the contrary to the symplectic case, the 
length of the Hamiltonian isotopy depends on the choice of a Hamiltonian function that generates the isotopy.
Finally, define the Hofer pseudo-norm on $\mathrm{Ham}(M,\pb)$ by
\[\nu(\Phi):=\inf\big\{l(f)\ :\ f\in C^{\infty}_{c}([0,1]\times M),\ \varphi_f^1=\Phi\big\}.\]
The following compatibility properties with the group structure are easily verified
\begin{itemize}
\item[\emph{(a)}] $\nu(\Phi) = \nu(\Phi^{-1})$,
\item[\emph{(b)}] $\nu(\Phi\circ\Psi)\leq \nu (\Phi)+\nu(\Psi)$,
\item[\emph{(c)}] $\nu(\Phi\circ\Psi\circ\Phi^{-1})=\nu(\Psi)$,
\end{itemize}
for $\Phi,\Psi\in \mathrm{Ham}(M,\pb)$ (see e.g.\ \cite[Theorem 1.1]{Rybicki}), but the non-degeneracy of $\nu$ is non-trivial: 
\begin{itemize}
\item[\emph{(d)}] $\nu(\Phi) = 0$ if and only if $\Phi=\mathrm{Id}$.
\end{itemize}

The main result of this article is the following.

\begin{theorem} \label{main-thm}
Let $(M,\pb)$ be a Poisson manifold. Then the Hofer pseudo-norm $\nu$ on $\mathrm{Ham}(M,\pb)$ is non-degenerate.
\end{theorem}

In the symplectic case, Hofer \cite{Hofer} proved non-degeneracy for the standard symplectic structure on $\mathbb{R}^{2n}$, then Polterovich \cite{Polterovich} extended it to a larger class of symplectic manifolds, and Lalonde and McDuff \cite{McDuff} completed the proof for all symplectic manifolds. All proofs rely on hard methods from symplectic topology. 

Below, we show that the Poisson case can be reduced to the sympletic case, by restricting to a symplectic leaf. This was first claimed by Sun and Zhang \cite{Sun}, in the setting of regular Poisson manifolds. Actually, in the proof they do not use regularity, but assume that the restriction of a compactly supported function to a leaf is compactly supported, however, without stating this explicitly. This property is equivalent to the leaves being closed submanifolds (which implies that they are embedded submanifolds, see e.g.\ \cite{FrMa}). This mistake was noticed by Rybicki \cite{Rybicki}, who obtained non-degeneracy for Poisson manifolds whose closed leaves form a dense set. Moreover, Rybicki \cite{Rybicki} proved non-degeneracy also for integrable Poisson manifolds, by using the displacement energy techniques on the symplectic groupoid. By adapting this proof to a sympelctic leaf, we obtain non-degeneracy in general. For the reader's convenience we will now briefly recall the definition and some properties of symplectic leaves of a Poisson manifold which are going to be used in the proof of Theorem \ref{main-thm}.

Consider the standard action of the Hamiltonian group $\mathrm{Ham}(M,\pb)$ on $M.$ The orbits of this action are called \emph{symplectic leaves}. Each symplectic leaf $L$ carries a unique smooth structure for which the inclusion $i: L \hookrightarrow M$ is an immersion, and moreover $L$ is an initial submanifold of $M.$\footnote{A submanifold $N$ of a manifold $M$ is called \emph{initial} if for every manifold $P$ and every smooth map $f:P \arr M$ such that $f(P) \subseteq N$ it holds that the induced map $f: P \arr N$ is smooth.} The tangent space at each point of $L$ is spanned by Hamiltonian vector fields. We equip $L$ with a canonical symplectic structure $\omega_L$ defined by
\begin{equation} \label{leaf_symp}
\omega_L(X_f, X_g) = \{f,g\}.
\end{equation}
Then the inclusion map $i: (L, \pb_L) \hookrightarrow (M, \pb)$ is a Poisson map, where $\pb_L$ is the Poisson structure induced by $\omega_L.$ For more details, we refer the reader to \cite[Section 1.3.4, p. 26]{Gengoux}.

We are now ready for the proof.

\begin{proof}[Proof of Theorem \ref{main-thm}]
Let $\Phi\in \mathrm{Ham}(M,\pb)$, $\Phi\neq \mathrm{Id}$, and fix $x_0\in M$ such that $\Phi(x_0)\neq x_0$. Let $i:L\to M$ be the symplectic leaf passing through $x_0$, where $i$ denoted the inclusion. 
Let $B\subset L$ be an open ball with compact closure such that $x_0 \in i(B)$ and 
\begin{equation} \label{displaced}
\Phi(i(B))\cap i(B)=\emptyset.
\end{equation}
Consider $f\in C^{\infty}_{c}([0,1]\times M)$ such that $\Phi=\varphi_f^1$. 
Since Hamiltonian isotopies preserve symplectic leaves, note that $\varphi^t_f(i(B)) \subseteq i(L)$. Choose a compactly supported smooth function $\lambda: L \to [0,1]$, such that $\lambda|_C= 1$, where $$C:=i^{-1}\Big(\cup_{t\in [0,1]}\varphi^t_f(i(B)) \Big)\subset L.$$ 
Define $g \in C_c^{\infty}([0,1] \times L)$ by 
$$g_t(y) := \lambda(y) \cdot f_t(i(y)).$$ 
Denote by $\varphi_g^t \in \mathrm{Ham}(L, \omega_L), t \in [0,1]$, the Hamiltonian flow of $g$, where $\omega_L$ is the symplectic form on $L$ (\ref{leaf_symp}). 
We need the following two claims to finish the proof.

\setcounter{theorem}{0}

\begin{claim} \label{claim_flows}
For every $t \in [0,1]$ we have that
$i \circ \varphi_g^t \vert_{B} = \varphi_f^t \circ i \vert_B.$
\end{claim}

\begin{proof}[Proof of Claim \ref{claim_flows}]
Fix $y\in B$. Since $i:L\hookrightarrow M$ is an initial submanifold, there is a (unique) smooth curve $\gamma:[0,1]\to L$ such that $i(\gamma(t))=\varphi_f^t(i(y))$. We need to show that $\gamma(t)=\varphi_g^t(y)$. Using that the inclusion $i: L \hookrightarrow M$ is a Poisson map it follows that:
\begin{equation}\label{Poisson map}
X^t_f\circ i(\gamma(t))= i_* X^t_{f \circ i}(\gamma(t)),\quad  \forall\ t \in [0,1].
\end{equation}
The fact that $\gamma(t)\in C$ and $\lambda|_C = 1$ implies that 
\begin{equation}\label{Other eq}
X^t_{f\circ i} (\gamma(t))=X^t_g (\gamma(t)) ,\quad  \forall\ t \in [0,1].
\end{equation}
Combining equations \eqref{Poisson map} and \eqref{Other eq}, we obtain that 
\[i_*(\frac{d}{dt} \gamma(t))=\frac{d}{dt}\varphi_f^t\circ i(y)=X^t_f\circ i(\gamma(t))=i_*(X^t_g(\gamma(t))),\quad  \forall\ t \in [0,1],\]
and since $i$ is an immersion, we obtain that $\gamma(t)$ is an integral line of $X_g^t$:
\[\frac{d}{dt} \gamma(t)=X^t_g(\gamma(t)),\quad  \forall\ t \in [0,1].\]
Since at $t=0$, $\gamma(0)=y$, Claim 1 follows.
\end{proof}

\begin{claim} \label{claim_lengths}
$l(g) \leq l(f).$
\end{claim}

\begin{proof}[Proof of Claim \ref{claim_lengths}]
We first consider the \textbf{case} when $M$ is compact.
By replacing $f_t$ by $f_t-f_t(x_0)$ we may assume that $f_t(x_0) = 0,$ for all $t \in [0,1].$ Note that the new function generate the same isotopy and that it has the same length as $f.$
Then we have that $\sup_{x\in M} f_t(x) \geq 0$ and therefore
\begin{equation} \label{sup}
\sup_{x\in M} f_t(x) \geq \sup_{y\in L}  \lambda(y) f_t(i(y)) = \sup_{y\in L} g_t(y), \quad t \in [0,1].
\end{equation}
Applying the same argument to $-f$, we get that
\begin{equation} \label{inf}
\inf_{x\in M} f_t(x) \leq \inf_{y\in L} g_t(y), \quad t \in [0,1].
\end{equation}
Combining inequalities (\ref{sup}) and (\ref{inf}) it follows that $l(g) \leq l(f)$ which
completes the proof of the claim in the case when $M$ is compact.

Now consider the \textbf{case} when $M$ is not compact. Then since $f$ is compactly supported we have that $\sup_{x\in M} f_t(x)  \geq 0$, for all $t \in [0,1].$ The same arguments as in the first case prove that the inequalities (\ref{sup}) and (\ref{inf}) hold also when $M$ is not compact. Hence
$l(g) \leq l(f)$ which completes the proof of Claim \ref{claim_lengths}.
\end{proof}

By (\ref{displaced}) and Claim \ref{claim_flows} it follows that $\varphi_g^1$ displaces $B \subseteq L.$ Hence $E(B)\leq l(g),$
where $E(B)$ is the displacement energy (see for example \cite[p. 469]{McDuffSalamon}) of the ball $B$ inside the symplectic manifold $(L, \omega_L).$
Then Claim \ref{claim_lengths} implies that $E(B)\leq l(f),$ and therefore
$$E(B)\leq \nu(\Phi).$$ 
Since $0<E(B)$ (see \cite[Theorem 1.1]{McDuff}) we conclude that $\nu(\Phi)>0.$ 
This completes the proof of Theorem \ref{main-thm}.
\end{proof}


\end{document}